# Peridynamic modeling for impact failure of wet concrete considering the influence of saturation

Liwei Wu[1], Dan Huang[1*], Qipeng Ma[1], Zhiyuan Li[1], Xuehao Yao[1]

1. College of Mechanics and Materials, Hohai University, Nanjing, 211100, China.

**Abstract**: In this paper, a modified intermediately homogenized peridynamic (IH-PD) model for analyzing impact failure of wet concrete has been presented under the configuration of ordinary state-based peridynamic theory. The meso-structural properties of concrete are linked to the macroscopic mechanical behavior in the IH-PD model, where the heterogeneity of concrete is taken into account, and the calculation cost does not increase. Simultaneously, the porosity of concrete is considered, which is implemented by deleting the bond between two material points, as well as the influence of porosity on the mechanical properties of concrete. Moreover, the effective bulk and shear modulus of cement mortar in wet concrete (saturated and unsaturated concrete) are calculated respectively. The dynamic model for wet concrete is described from three aspects: strength, dynamic increase factor, and equation of state. Validation of the proposed model is established through analyzing some benchmark tests and comparing with the corresponding experiment and other available numerical results.

## 1. Introduction

As one of the most widely used engineering materials, concrete has attracted the focus of concern of researchers from various engineering fields for a long time. In many engineering structures, concrete materials work in the water environment, for example, in the dam, pier, aqueduct, various offshore buildings, etc. With the effect of external water pressure, the water content in concrete will be different, which will accordingly take effects on the behavior and properties of concrete materials and structures. This brings new challenging but important issues to the studies on the mechanical behavior of concrete materials and structures, especially when dynamic and failure problems are concerned.

With the rapid development of computer technology, numerical methods have become a significant tool for solving various mechanical problems. Among the existing numerical methods, most of them are based on the classical continuum mechanics theory. Nevertheless, the continuity hypothesis in the continuum mechanics theory is contradictory to reality, and the control equation cannot be directly used to the discontinuous problems. Some remedy techniques and numerical approaches were developed (Cundall and Strack, 2008; Ma et al., 2010; Rizzo and Frank, 1967), but the fundamental contradiction between the discontinuities and continuity hypothesis is not solved. The solution for the dynamic crack propagation and impact failure is still challenging. Different from other traditional numerical methods, a non-local theory, named peridynamics (PD) (Silling, 2000), was proposed, which owns promising advantages in analyzing the discontinuous problems (Bobaru and Zhang, 2015; Ha and Bobaru, 2010). Peridynamics, presented as a reformulation of the traditional continuum mechanics theory, adopts the integral-differential equation of motion instead of the partial-derivatives based one on conventional theories. Meanwhile, peridynamics includes the description of damage and fracture (Foster et al., 2011; Silling and Askari, 2005), which means that the crack could propagate spontaneously without any additional fracture criteria. Such a model has the natural theoretical advantages on analyzing dynamic crack growth (Gu et al., 2016; Wang et al., 2018; Zhou et al., 2021), impact (Ren et al., 2017; Wang et al., 2019; L Wu, Xu, et al., 2020; Wu et al., 2021) and blast (Fan et al., 2016; Fan and Li, 2017; Zhu and Zhao, 2021) failure, and other discontinuous problems with non-local effects (He et al., 2021; Li et al., 2021; Ni et al., 2020; Wang et al., 2016, 2017; Zhang and Qiao, 2020).

Many studies on the mechanical behavior of concrete with peridynamics have been conducted in recent years. Yaghoobi et al. (Yaghoobi and Chorzepa, 2015) predicted the response of fiber reinforced concrete, using a semi-discrete peridynamic method to model the fiber reinforcement. Yang et al. (Yang et al., 2018) proposed a new damage model to quantitatively investigate the mode-I crack propagation in concrete. Huang et al. (Huang et al., 2015) presented an extended bond-based peridynamic approach for the quasi-static mechanical behavior of concrete materials and structures. Xu et al. (Xu et al., 2021) simulated the damage processes of prefabricated beam post-cast with steel fiber reinforced high-strength concrete based on the micro-polar peridynamics. Zhang et al. (Zhang et al., 2021) simulated the relative slip and interface damages between the concrete and steel in RC structures by a novel coupled axial-shear interaction (ASI) bond-slip model. Gu et al. (Gu et al., 2016) studied the elastic wave dispersion, propagation, and fracture of the concrete through a split Hopkinson pressure bar model. In our previous work, the impact failure under high-velocity impact conditions was analyzed (Wu et al., 2019; L Wu, Huang, et al., 2020).

The foregoing studies on concrete are mostly performed with the homogenized model in peridynamics, ignoring the meso-structural characteristics of concrete, which may cause inaccurate descriptions of the heterogeneities-driven fracture and failure(P Wu, Zhao, et al., 2020). Several mesoscopic models for concrete have been proposed and verified (Dong et al., 2021; Li and Guo, 2018; Shi et al., 2021). Nevertheless, these models are restricted due to their huge computation cost, especially for three-dimensional dynamic problems. Recently, the IH-PD model (Chen et al., 2019) was proposed for porous materials (Chen et al., 2019) and functionally graded materials (Mehrmashhadi et al., 2019), and it also has been successfully applied to concrete materials (P Wu, Zhao, et al., 2020; Wu and Huang, 2022; Zhao et al., 2020).

In the present work, the impact failure of wet concrete is modeled and simulated based on the IH-PD model. The concrete porosity is taken into account, which is implemented by deleting the bond between two material points. The influence porosity has on the mechanical behavior of concrete is considered as well. Besides, the relationship between the water content and the dynamic behavior in wet concrete is established. On the one hand, the effective mechanical properties of cement mortar in wet concrete are considered, according to the equivalent method with the two-phase spherical model. On the other hand, the dynamic model of wet concrete with different saturations is described from three aspects: strength, dynamic increase factor, and equation of state. We verify the applicability of the presented model and approach on accurately characterizing the impact failure of wet concrete.

The paper is organized as follows: In Section 2, the ordinary state-based peridynamic theory is briefly introduced. The modified IH-PD model of wet concrete is introduced in Section 3. In Section 4, validation of the proposed model and approach is discussed through several typical numerical examples. Finally, concluding remarks are summarized in Section 5.

## 2. Methodology

### 2.1. Brief review of ordinary state-based peridynamics

In the state-based peridynamic theory (Silling et al., 2007), the motion equation is given as:

$$\rho \ddot{\mathbf{u}}(\mathbf{x},t) = \int_{H_\mathbf{x}} \left\{ \underline{\mathbf{T}}[\mathbf{x},t]\langle \mathbf{x}' - \mathbf{x} \rangle - \underline{\mathbf{T}}[\mathbf{x}',t]\langle \mathbf{x} - \mathbf{x}' \rangle \right\} dV_{\mathbf{x}'} + \mathbf{b}(\mathbf{x},t). \qquad (1)$$

where $\rho$ and $\ddot{\mathbf{u}}$ is the material mass density and the acceleration, $dV_{\mathbf{x}'}$ is the volume for $\mathbf{x}'$, and $\mathbf{b}$ denotes the external body force density. $\underline{\mathbf{T}}[\mathbf{x},t]$ and $\underline{\mathbf{T}}[\mathbf{x}',t]$ are force vector states that represent the relationship between material points at time $t$. The integral domain $H_{\mathbf{x}}$ is the horizon region, denoted as $H_{\mathbf{x}} = \{0 < |\mathbf{x}-\mathbf{x}'| < \delta\}$ (as shown in Fig. 1).

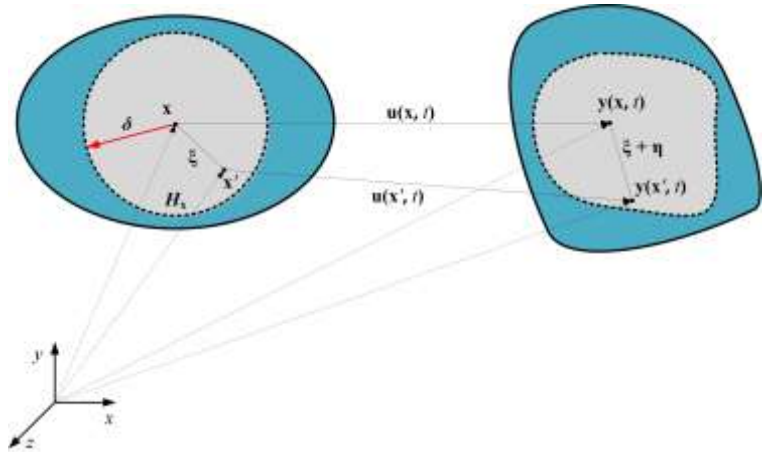

Fig. 1. The interaction between material points in the peridynamic theory.

In the peridynamic theory, there are three types of typical models, namely bond-based, ordinary state-based, and non-ordinary peridynamic models, respectively. The difference is the direction of bond force state, as shown in Fig. 2. Besides, the non-ordinary peridynamic model has inherent stability issues caused by zero-energy modes and could result in some unrealistic errors in numerical results. Therefore, the ordinary peridynamic model is used in this study.

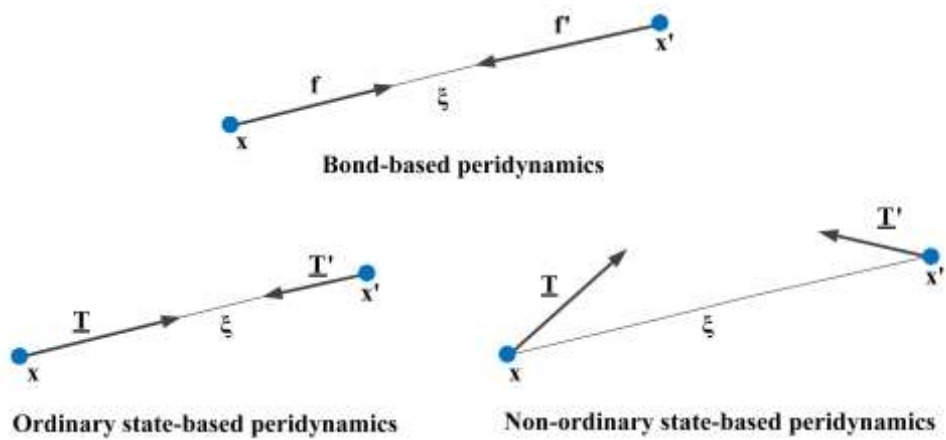

Fig. 2. Three types of typical peridynamic models.

For ordinary state-based peridynamics, the force vector state $\underline{\mathbf{T}}$ satisfies:

$$\underline{\mathbf{T}} = \underline{t}\,\underline{\mathbf{M}}, \tag{2}$$

where $\underline{t}$, $\underline{\mathbf{M}}$ denotes a scalar force state and a unit vector in the deformed configuration, respectively.

For elastic material, the strain energy density $W$ is introduced as (Silling et al., 2007):

$$W\left(\theta, \underline{e}^d\right) = \frac{k\theta^2}{2} + \frac{\alpha}{2}\left(\underline{\omega}\,\underline{e}^d\right) \cdot \underline{e}^d, \tag{3}$$

where $k$ and $\alpha$ are positive constants.

According to Eq. (3), the scalar force state in Eq. (2) is denoted as:

$$\underline{t} = \frac{3k\theta}{m}\underline{\omega}\underline{x} + \alpha\underline{\omega}\underline{e}^d, \tag{4}$$

$$\theta(\underline{e}) = \frac{3}{m}(\underline{\omega}\underline{x}) \cdot \underline{e}, \tag{5}$$

$$m = (\underline{\omega}\underline{x}) \cdot \underline{x}, \tag{6}$$

where $\underline{e}$ is the extension scalar state, denoted by:

$$\underline{e} = \underline{y} - \underline{x}, \quad \underline{y} = |\underline{\mathbf{Y}}|, \quad \underline{x} = |\underline{\mathbf{X}}|, \tag{7}$$

where $\underline{x}$ and $\underline{y}$ is the scalar state in reference and deformation position.

## 2.2. The damage and failure criterion

In peridynamics, once the stretch of a bond exceeds a critical value, that bond is regarded as being broken.

A boolean function is introduced here:

$$\mu(t,\xi) = \begin{cases} 1 & s_0^{fc} < s(t',\xi) < s_0^{ft},\ 0 < t' < t \\ 0 & else \end{cases}, \tag{8}$$

where $s$ denotes the bond stretch, $s_0^{fc}$ and $s_0^{ft}$ is a compressive and tensile critical stretch ($s_0^{ft} = \sqrt{\dfrac{5G_0}{9k\delta}}$ (Silling and Askari, 2005), $s_0^{fc} = \dfrac{\sigma_c}{E}$ (L Wu, Huang, et al., 2020)).

The bond is intact only if $\mu = 1$. Fig. 3 shows the relationship between the bond force and bond stretch.

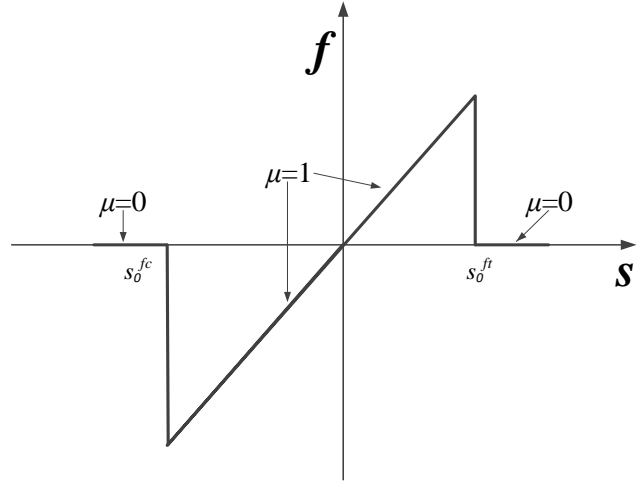

Fig. 3. Relationship between the bond force and bond stretch.

In peridynamic theory, the damage value $D$ of a material point is generally denoted by the ratio of number of broken bonds to number of total bonds:

$$D(x,t) = 1 - \dfrac{\int_{H_x} \mu(\xi,t) dV_{x'}}{\int_{H_x} dV_{x'}} \ . \tag{9}$$

If the value of $D$ is zero, corresponds to the intact material. If the value of $D$ is one, corresponds to the fully damaged material.

### 2.3. Contact model

As shown in Fig. 4, the contact force is independent of the original distance between material points, only depending on the current relative position (Macek and Silling, 2007).

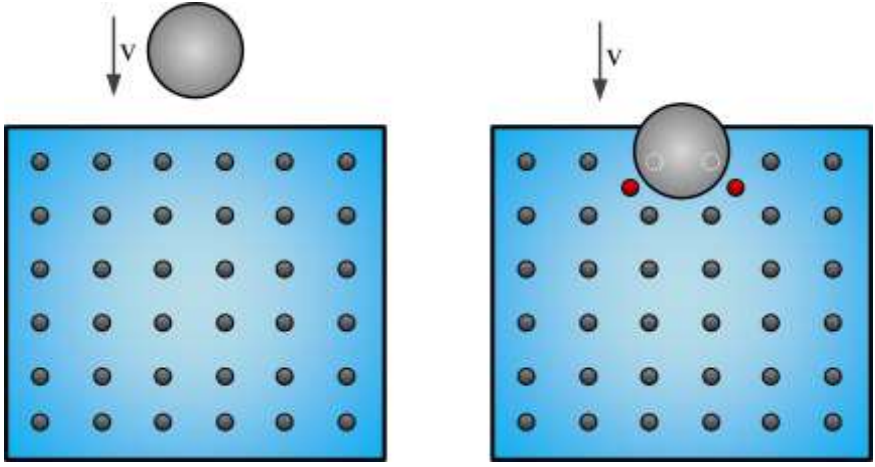

Fig. 4. Impact of material points with the contact model.

The contact force is:

$$f_s(\mathbf{y}_p, \mathbf{y}_i) = \min\{0, \frac{c_{sh}}{\delta}(\|\mathbf{y}_p - \mathbf{y}_i\| - d_{pi})\} \frac{\mathbf{y}_p - \mathbf{y}_i}{\|\mathbf{y}_p - \mathbf{y}_i\|}, \qquad (10)$$

where $\mathbf{y}$ denotes the position of material point. $c_{sh} = 15c$, and $c$ is a peridynamic constant corresponding to the classical continuum mechanics bulk modulus $k$, which is found by equating the strain energy under isotropic extension from continuum mechanics to the energy density (Macek and Silling, 2007). $d_{pi}$ denotes the critical distance:

$$d_{pi} = \min\{0.9\|\mathbf{x}' - \mathbf{x}\|, 1.35|\Delta x|\}. \qquad (11)$$

Considering the contact force, Eq. (1) could be formulated as:

$$\rho \ddot{\mathbf{u}}(\mathbf{x},t) = \int_{H_x} \{\underline{\mathbf{T}}[\mathbf{x},t]\langle \mathbf{x}' - \mathbf{x}\rangle - \underline{\mathbf{T}}[\mathbf{x}',t]\langle \mathbf{x} - \mathbf{x}'\rangle\} dV_{\mathbf{x}'} + (\mathbf{b}(\mathbf{x},t) + \mathbf{f}_s(\mathbf{x},t)). \qquad (12)$$

This contact model has reasonable results in elastic conditions, but is not suitable for plastic impact problems, where energy dissipation caused by plastic deformation is significant (Wu and Huang, 2021).

## 3. Modified IH-PD model for wet concrete

### 3.1 IH-PD model

In an IH-PD model, the mechanical characteristics of concrete meso-structure are linked to the macroscopic mechanical behavior, implemented through the bond level. This model does not describe the actual structure in concrete, only aiming to maintain the heterogeneity information about the meso-structure and reduce the computational cost. Here the concrete is regarded as a three-phase model, consisting of cement mortar, aggregates, and ITZ. As shown in Fig. 5, there are six types of combined bonds. When two material points locate in the same phase, A-A bond, C-C bond, and I-I bond are formed. When two material points locate in a different phase, A-C bond, A-I bond, and C-I bond are formed. Besides, the mechanical characteristic of A-C bond is assumed to be equivalent to the homogenous concrete model in this study, while A-I and C-I bonds are regarded to be identical to I-I bond.

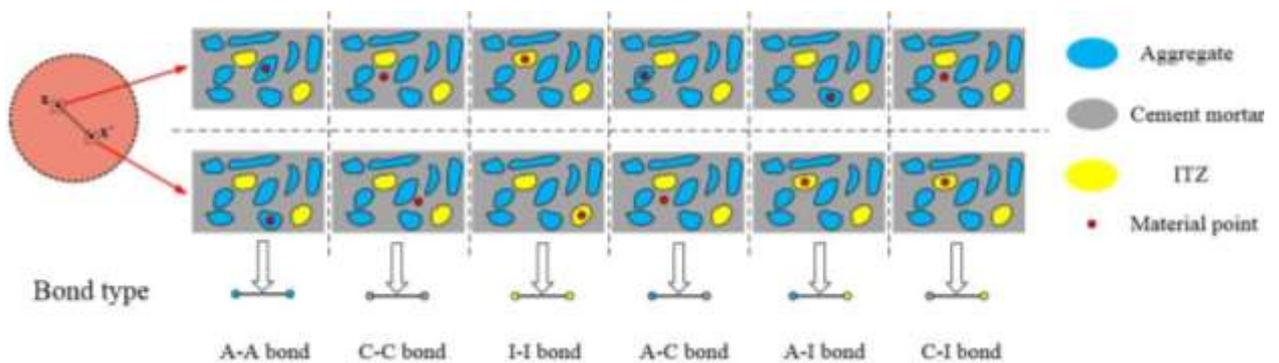

Fig. 5. Description of different types of the combined bond.

To have a better understanding of the real concrete material, the porosity is represented in this study, which can be expressed by the pre-broken bonds between material points. Similar to the foregoing process, we first determine the type of bond, then randomly determine whether the bond is broken or not. The pore bond does not need to be deleted, just simply broken in the pre-processing process, as shown in Fig. 6.

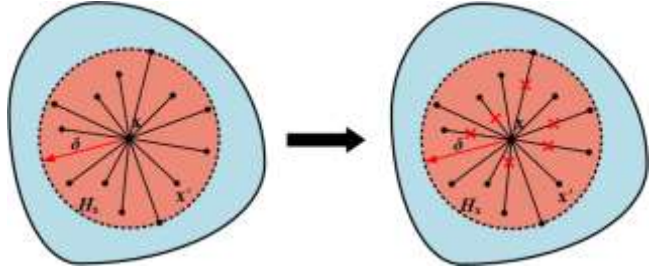

Fig. 6. The bond broken procedure in the IH-PD model.

To perform the number of pre-broken bonds for arbitrary material point, a pre-damage index $d$ is defined as:

$$d(\mathbf{x}) = \frac{\varphi(\mathbf{x})}{\varphi_c}, \tag{13}$$

where $\varphi(\mathbf{x})$ is the concrete porosity, and $\varphi_c$ is the critical concrete porosity. When $\varphi(\mathbf{x})$ reaches $\varphi_c$ (assumed to be 1 in this study), all bonds associated with that material point should be broken. The pre-processing algorithm for implementing pre-broken bonds is shown in Fig. 7. For example, the volume fraction is chosen to be 40% (aggregate), 55% (cement mortar), and 5% (ITZ), concrete porosity is 10%. The proportion for the aggregate bond is 16%, cement mortar bond is 74.25%, ITZ bond is 9.75%, and pre-broken bond is 19%.

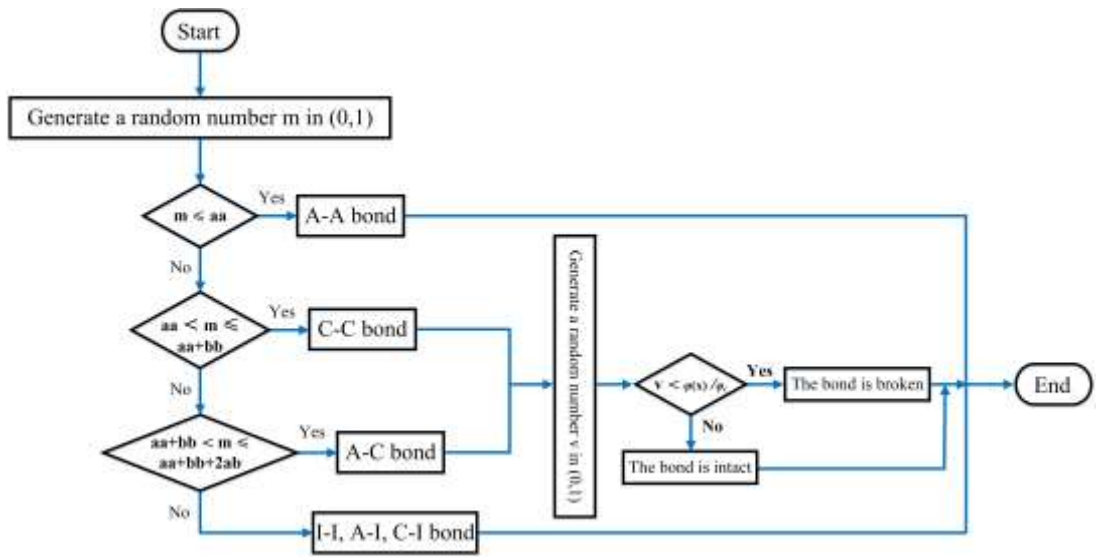

Fig. 7. The flowchart of a pre-processing algorithm for the IH-PD model considering concrete pores.

Due to the existence of internal pores, the macroscopic mechanical properties of concrete would vary accordingly. In this study, referring to the method in literature (Jin et al., 2012), cement mortar with pores is regarded as a two-phase medium, composed of mortar matrix (zero porosity) and pores. It is assumed that the porosity is the volume fraction of all pores in cement mortar, and it remains to be constant. The effective bulk modulus and shear modulus of cement mortar with pores are given by (Jin et al., 2012):

$$K^* = \frac{4K_m \mu_m (1-\varphi)}{(4\mu_m + 3K_m \varphi)},$$
$$\mu^* = \mu_m (1-\varphi^2)$$
(14)

where $K_m$, $\mu_m$ is the bulk modulus and shear modulus of mortar matrix, $\varphi$ is the given material porosity.

### 3.2 Effective mechanical properties of wet concrete

### 3.2.1 Saturated concrete

For saturated concrete, the cement mortar is represented to be a two-phase model, composed of mortar matrix and free water (Fig. 8(a)). Since pores are filled with water, the volume fraction of pore water is equal to the porosity, and the simplified physical model is shown in Fig. 8(b). The radius of the inner and outer spheres in Fig. 8(b) is $a$ and $b$, and the porosity can be expressed as $\varphi = a^3 / b^3$.

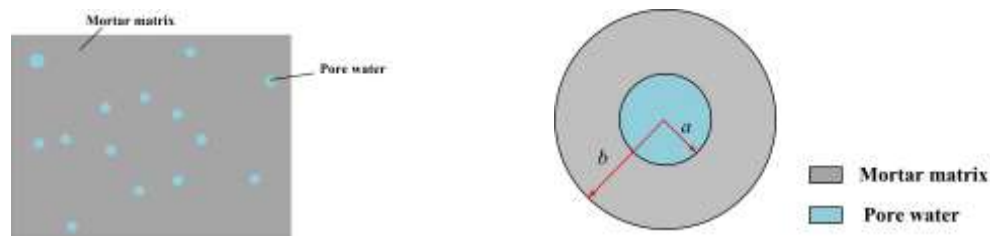

(a) Schematic diagram of cement mortar containing pore water.

(b) The simplified two-phase sphere model.

Fig. 8. The simplified model of cement mortar in saturated concrete.

To establish the relationship between pore water and porosity, an equivalent analysis is made here, as shown in Fig. 9. The bulk modulus of pore water is $K_w$ ($K_w < K_m$). The pore water can be equivalent to a two-phase model (Jin et al., 2012), indicating that $K_w$ is identical to the bulk modulus of concrete with porosity $\varphi_1$. According to Eq. (14), $K_w$ satisfies the following relationship:

$$K_w = \frac{4K_m \mu_m (1-\varphi_1)}{4\mu_m + 3K_m \varphi_1}. \tag{15}$$

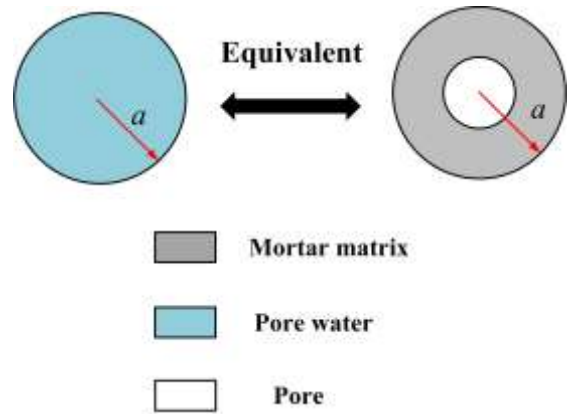

Fig. 9. The equivalent model of pore water and mortar matrix with a pore.

The porosity can be calculated from Eq. (15) as:

$$\varphi_1 = \frac{4\mu_m (4K_m - K_w)}{K_m (3K_w + 4\mu_m)}. \tag{16}$$

The pore water in cement mortar (shown in Fig. 8(b)) could be replaced with the two-phase model for mortar matrix, and the final two-phase model is established (shown in Fig. 10). The total porosity $\varphi_2$ for the final two-phase model can be evaluated as:

$$\varphi_2 = \frac{a^3 \cdot \varphi_1}{b^3} = \varphi \varphi_1. \tag{17}$$

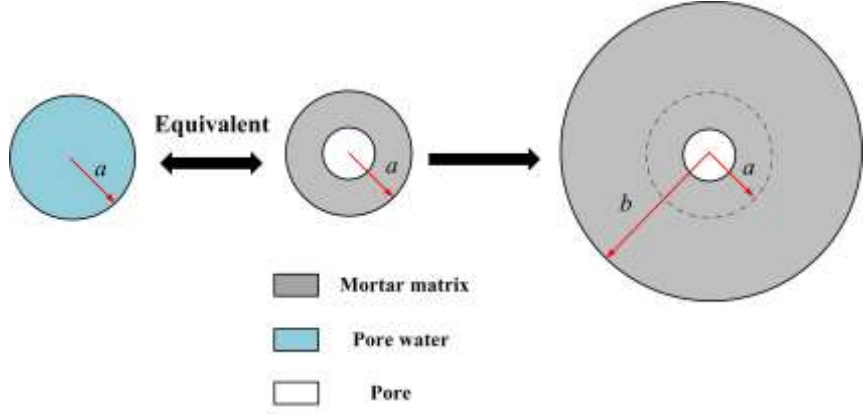

Fig. 10 Complete equivalence of the saturated concrete.

Substituting Eq. (17) into Eq. (16), the effective bulk modulus in saturated concrete, depicted in Fig. 8(b), is calculated by:

$$K_w^* = \frac{4K_m\mu_m(1-\varphi_2)}{4\mu_m+3K_m\varphi_2} = \frac{4K_m\mu_m(1-\varphi\varphi_1)}{4\mu_m+3K_m\varphi\varphi_1}. \tag{18}$$

Generally speaking, the shear resistance of water could be negligible. However, when subjected to dynamic load, the viscosity of pore water is a significant factor in the mechanical properties of saturated concrete, except for the inertia effect. The influence viscosity of pore water has on the shear modulus of saturated concrete cannot be directly ignored. According to (Wang and Li, 2007), the relationship between effective shear modulus and porosity in saturated concrete is established by:

$$\mu_w^* = (1+f_1\varphi^2+f_2\varphi)(1-\varphi^2)\mu_m. \tag{19}$$

### 3.2.2 Unsaturated concrete

For the wet concrete in a water environment, it takes a long time to reach the saturated state. For example, Chatterji (Chatterji, 2004) measured that the concrete structure was still dry even had been stored in the wet environment for 222 days. Persson (Persson, 1997) tested the saturation of concrete was 98%, not fully saturated, even the concrete samples had been cured in water for 450 days. It is extremely meaningful to analyze the mechanical characteristics of unsaturated concrete.

The concrete saturation refers to the volume ratio of the pore water, but this definition does not apply to unsaturated concrete, where pores are not just filled with water, as shown in Fig. 11. There are three types of pore in unsaturated concrete, satisfying the following relationship:

$$V_p = V_{p-sat} + V_{p-unsat} + V_{p-dry}, \tag{20}$$

where $V_{p-sat}$, $V_{p-unsat}$, and $V_{p-dry}$ is the volume of saturated pores, unsaturated pores, and dry pores, respectively. $V_p$ is the total volume of wet concrete.

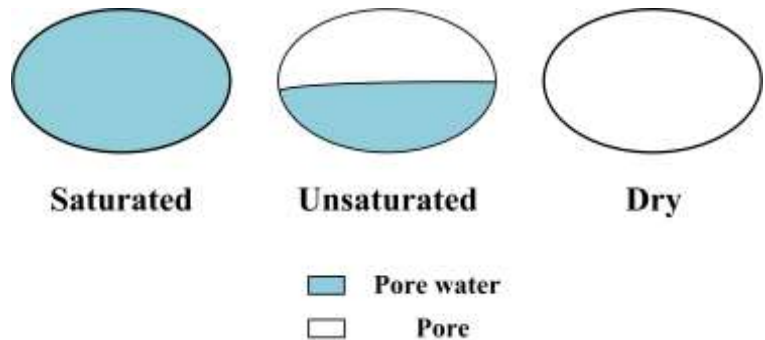

Fig. 11. The type of pore in saturated, unsaturated, and dry concrete.

The corresponding total porosity $\varphi$ is:

$$\varphi = \varphi_{sat} + \varphi_{unsat} + \varphi_{dry}, \tag{21}$$

where $\varphi_{sat}$ is saturated porosity, $\varphi_{unsat}$ is unsaturated porosity, and $\varphi_{dry}$ is dry porosity.

In this section, the cement mortar is composed of a mortar matrix, free water, and pores. To obtain the effective modulus of it, just considering porosity is not enough. The influence of saturation on unsaturated concrete should be considered.

Firstly, the two-phase model of dry pores, unsaturated pores, and mortar matrix is established, as shown in Fig.12. According to the last section, the $K_1^*$ and $\mu_1^*$ of equivalent body A are calculated by:

$$K_1^* = \frac{4K_m\mu_m(1-\varphi_1')}{4\mu_m + 3K_m\varphi_1'},$$
$$\mu_1^* = \mu_m(1-\varphi_1'^2)$$
(22)

where $\varphi_1'$ is the porosity of equivalent body A, described by:

$$\varphi_1' = \frac{V_{dry} + \frac{1}{2}V_{unsat}}{V - V_{sat} - \frac{1}{2}V_{unsat}} = \frac{\varphi_{dry} + \frac{1}{2}\varphi_{unsat}}{1 - \varphi_{sat} - \frac{1}{2}\varphi_{unsat}} = \frac{\varphi - \varphi \cdot w}{1 - \varphi \cdot w},$$
(23)

where $w$ is the saturation of unsaturated concrete. Since the water content in unsaturated pores is not easy to be determined, it is assumed to be equal to 1/2 in the numerical simulation, according to the average hypothesis. Therefore, the volume and porosity of unsaturated pores in Eq. (23) are multiplied by 1/2.

Substituting Eq. (23) into Eq. (22):

$$K_1^* = \frac{4K_m\mu_m\left(1 - \left(\frac{1-\varphi}{1-\varphi \cdot w}\right)\right)}{4\mu_m + 3K_m \cdot \left(\frac{\varphi - \varphi \cdot w}{1-\varphi \cdot w}\right)}.$$
$$\mu_1^* = \mu_m\left(1 - \left(\frac{\varphi - \varphi \cdot w}{1-\varphi \cdot w}\right)^2\right)$$
(24)

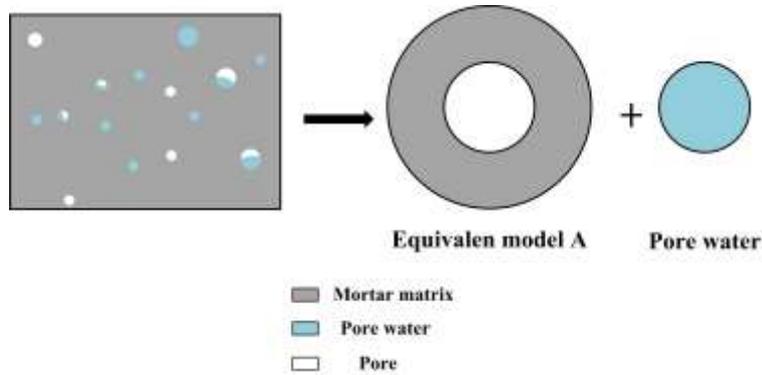

Equivalen model A    Pore water

☐ Mortar matrix
☐ Pore water
☐ Pore

Fig. 12 The equivalent model for unsaturated concrete.

Secondly, the two-phase model of pore water and equivalent body A is established,

and it could be regarded as an equivalent model of saturated concrete. The following procedure is similar to the last section. Substituting Eq. (24) into Eq. (18) and Eq. (19), the effective modulus of unsaturated concrete can be calculated as:

$$K^*_{unsat} = \frac{4K_1^* \mu_1^* (1-\varphi_2)}{4\mu_1^* + 3K_1^* \varphi_2} = \frac{4K_1^* \mu_1^* \left(1 - \left(\varphi_{sat} + \frac{1}{2}\varphi_{unsat}\right)\varphi_1\right)}{4\mu_1^* + 3K_1^* \left(\varphi_{sat} + \frac{1}{2}\varphi_{unsat}\right)\varphi_1} = \frac{4K_1^* \mu_1^* (1 - \varphi \cdot w \cdot \varphi_1)}{4\mu_1^* + 3K_1^* \varphi \cdot w \cdot \varphi_1}$$

$$\mu^*_{unsat} = \mu_1^* \left(1 - \left(\varphi_{sat} + \frac{1}{2}\varphi_{unsat}\right)^2\right) = \mu_1^* \left(1 - (\varphi \cdot w)^2\right)$$

(25)

### 3.3 Dynamic model for wet concrete

Experimental tests have shown that saturation has a considerable effect on quasi-static and dynamic behavior of wet concrete. Meanwhile, it is found that the unconfined strength of wet concrete is slightly lower than that of dry concrete when subjected to quasi-static loading conditions. Through the experimental measurements, Zhao et al. (Zhao and Wen, 2018) presented an empirical equation of the quasi-static strength of wet concrete, and the relationship between the strength of wet concrete and its saturation can be expressed as follows:

$$\begin{aligned} \frac{f_c^w}{f_c} &= -0.2w + 1.0 \\ \frac{f_t^w}{f_t} &= -0.2w + 1.0 \end{aligned},$$

(26)

where $f_c$, $f_t$ are quasi-static compressive and tensile strength of dry concrete. $f_c^w$, $f_t^w$ are quasi-static strength of wet concrete. $w$ is the saturation of unsaturated concrete.

According to Eq. (8), the relationship between the critical stretch for wet concrete and concrete strength is presented as:

$$s_0^{tw} = f_t^w / E_w, \qquad s \geq 0$$
$$s_0^{cw} = f_c^w / E_w, \qquad s < 0 \tag{27}$$

where $s_0^{cw}$, $s_0^{tw}$ are the compressive and tensile critical stretch value, $E_w$ is Young's modulus of wet concrete.

The experimental tests (Malecot et al., 2019; Vu et al., 2009) demonstrated that wet concrete has more sensitive characteristics to strain rate than dry concrete. Thus it is extremely significant to consider its rate effect when subjected to dynamic loadings. For dry concrete, the rate effect under tension and compression conditions is described by

$$DIF_c^{dry} = \left[1 + C \ln(\dot{s}_c)\right]$$

$$DIF_t^{dry} = \begin{cases} \left(\dfrac{\dot{s}_t}{\dot{\varepsilon}_t}\right)^{\zeta}, & \dot{s}_t \leq 30s^{-1} \\ \beta\left(\dfrac{\dot{s}_t}{\dot{\varepsilon}_t}\right)^{\frac{1}{3}}, & \dot{s}_t > 30s^{-1} \end{cases}, \tag{28}$$

where $\dot{s}_c$, $\dot{s}_t$ are compressive and tensile strain rates, and $C$, $\zeta$ are corresponding positive material constants. To depict the influence of saturation on the dynamic behavior of wet concrete, Zhao et al. (Zhao and Wen, 2018) presented an empirical formula to define the ratio of the wet concrete dynamic increase factor to that of dry concrete:

$$g_t(\dot{s}, w) = \frac{DIF_t^w}{DIF_t^{dry}} = \begin{cases} 1, & \dot{s} \leq 10^{-5} \\ -A^{-\lg \dot{s} - 5} + 2, & \dot{s} > 10^{-5} \end{cases}, \tag{29}$$

$$g_c(\dot{s}, w) = \frac{DIF_c^w}{DIF_c^{dry}} = \begin{cases} 1, & \dot{s} \leq 10^{-5} \\ \dfrac{DIF_t^w \cdot f_t^{dry} / f_c^{dry} + 1}{DIF_t^{dry} \cdot f_t^{dry} / f_c^{dry} + 1}, & \dot{s} > 10^{-5} \end{cases}, \tag{30}$$

where $A$ is a material parameter, and it satisfies:

$$A = 1 + Cw, \tag{31}$$

where the value $C$ is 0.15.

According to Eq. (29) and Eq. (30), the dynamic increase factor of wet concrete can be denoted as:

$$DIF_c^w = \left[1 + C\ln(\dot{s}_c)\right] \cdot g_c(\dot{s}, w)$$

$$DIF_t^w = \begin{cases} \left(\dfrac{\dot{s}_t}{\dot{\varepsilon}_t}\right)^\zeta, & \dot{s} \leq 10^{-5} \\ \left(\dfrac{\dot{s}_t}{\dot{\varepsilon}_t}\right)^\zeta \cdot g_t(\dot{s}, w), & 10^{-5} < \dot{s}_t \leq 30 s^{-1} \\ \beta\left(\dfrac{\dot{s}_t}{\dot{\varepsilon}_t}\right)^{\frac{1}{3}} \cdot g_t(\dot{s}, w), & \dot{s}_t > 30 s^{-1} \end{cases}. \tag{32}$$

In the wet concrete, pores are full of free water, and the pore water is assumed to be in undrained condition when subjected to impacting load. The free water cannot flow from pores, resulting in the corresponding pores not collapsing. Therefore, the consolidation pressure of wet concrete $P_{lock}^w$ is lower than dry concrete $P_{lock}$, as shown in Fig. 13.

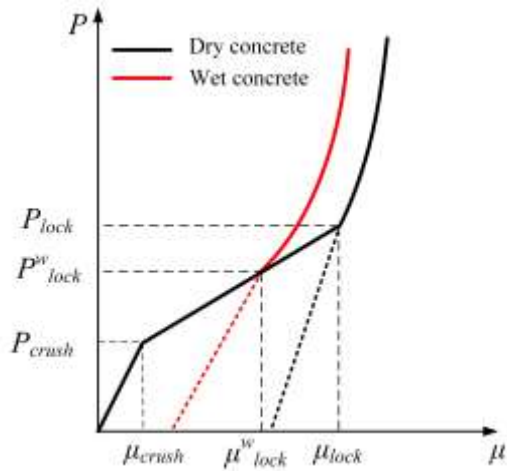

Fig. 13. The equation of state for dry concrete and wet concrete.

Because free water can stay in pores and bear partial hydrostatic pressure, the hydrostatic pressure of wet concrete $P'$ is composed of solid part $P_{dry}$ and free water part $P_w$ (Huang et al., 2020):

$$P' = P_{dry} + bP_w, \tag{33}$$

$$P_{dry} = k_1 \bar{\mu} + k_2 \bar{\mu}^2 + k_3 \bar{\mu}^3, \tag{34}$$

$$P_w = \frac{\rho_0 L^2 \mu' \left(1 + \left(1 - \frac{\gamma_0}{2}\right)\mu' - \frac{\alpha}{2}\mu'^2\right)}{1 - (S_1 - 1)\mu' - S_2 \frac{\mu'^2}{\mu'+1} - S_3 \frac{\mu'^3}{(\mu'+1)^2}} + (\gamma_0 + \alpha\mu')E, \tag{35}$$

where $b$ is the Biot coefficient, $b = 1 - (1-\varphi)^3$. $\bar{\mu} = \frac{\mu - \mu_{crush}}{\mu_{lock} - \mu_{crush}}$, $\mu$ is the volumetric strain. $\mu_{crush}$, $\mu_{lock}$, and $\mu'_{lock}$ are the volumetric strain under $P_{crush}$, $P_{lock}$, and $P_{lock}^w$. $k_1$, $k_2$, and $k_3$ are material parameters of concrete. Moreover, other constants are $\rho_0 = 1000 \text{kg/m}^3$, $L = 1480 \text{m/s}$, $S_1 = 2.56$, $S_2 = 1.986$, $S_3 = 1.2268$, $\gamma_0 = 0.35$, $E = 1.89 \times 10^6 \text{J/m}^3$, and $\alpha = 0$ (Huang et al., 2020).

## 4. Numerical verification and demonstration

### 4.1 Effective Young's modulus in a wet concrete target

The effective Young's modulus of a three-dimensional concrete slab is simulated here, where the elastic wave is generated by suddenly applying a pulse force on the concrete slab. The geometric information is shown in Fig. 14. The force pulse is 1MPa and suddenly drops to 0 at 5μs. Young's modulus of each phase, aggregate, cement, and ITZ are 70GPa, 15Gpa, and 7.5GPa, respectively (P Wu, Zhao, et al., 2020). Material properties of each phase are Besides, the volume fraction of aggregate and cement mortar is 40% and 57%, and this will be used in this paper unless otherwise specified. For peridynamic modeling, the horizon is $\delta=40\text{mm}$, and $\delta=4\Delta\text{x}$. The total number of nodes in this case is 61,206.

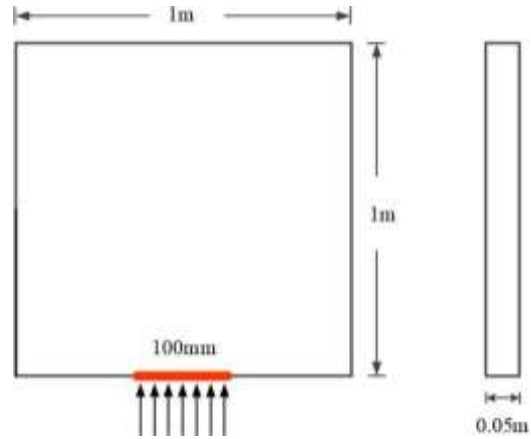

Fig. 14. The geometric information condition in a wet concrete slab.

The effective Young's modulus with different porosities, calculated according to wave velocity $C_w$ (Lee et al., 2017), is shown in Fig. 15, where $E_m$ is Young's modulus when the porosity is 0. In this study, the wave velocity is an average velocity in the vertical direction, which is determined by tracking the displacement of the appointed point over a specific time step. The porosity is 0, 10 %, 30 %, 50 %, and 70 %, respectively. In Fig. 15, it can be observed that the porosity increment dramatically softens Young's modulus. The proposed model provides an accurate prediction of the effective Young's modulus with different porosities, which appears to be experimentally repeatable (Wu, 2012). The good agreement indicates that the presented model could accurately describe the quantitative relationship between the effective Young's modulus and porosity for wet concrete.

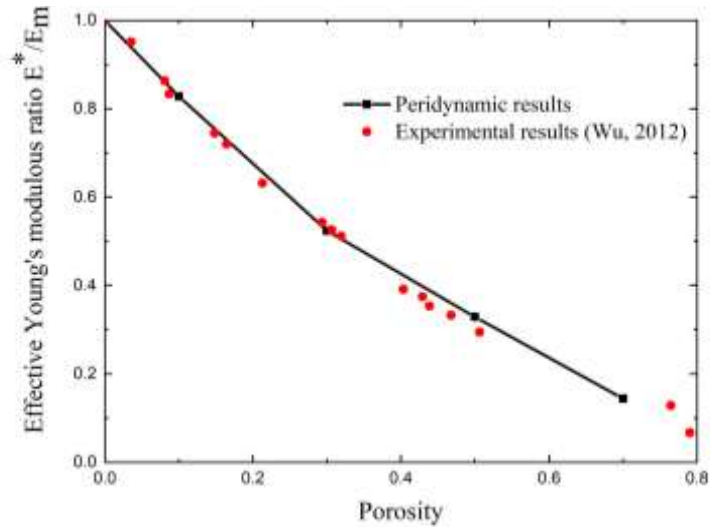

Fig. 15. Effective Young's modulus obtained by experimental measurements and peridynamic predictions with different porosities.

Fig. 16 illustrates the effective Young's modulus with different porosities, while concrete is fully saturated in Fig. 16(a), and the other one is completely dry in Fig. 16(b). The IH-PD results are highly consistent with the corresponding experimental data (Yaman et al., 2002) and numerical results from other models (Jin et al., 2012; Wang and Li, 2007), demonstrating that the proposed model has good reliability and accuracy for simulating the effective Young's modulus of saturated and dry concrete. Besides, the effective Young's modulus of saturated and dry concrete decreases dramatically when the concrete porosity increases, but the latter decreases faster than the former (same value of starting point, different value of ending point). The reason is that the pore water in saturated concrete would restrict the deformation of the surrounding mortar matrix, leading to an improvement of the concrete stiffness.

To study the quantitative relationship between effective Young's modulus and porosity in wet concrete, the numerical results for dry, unsaturated, and saturated concrete are depicted in Fig. 17, respectively. There are five different saturations. The effective Young's modulus gradually decreases when the porosity increases, no matter what is the saturation. For the same porosity, the effective Young's modulus of wet concrete is greater when the saturation is higher, and the value of saturated concrete is the maximum among different saturations.

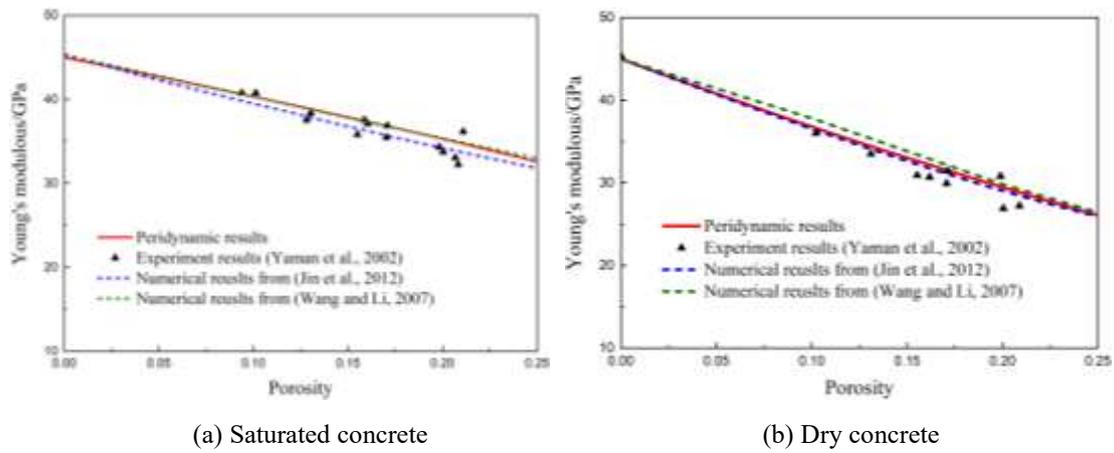

(a) Saturated concrete  (b) Dry concrete

Fig. 16. Comparison of the effective Young's modulus with different porosities in saturated and dry concrete.

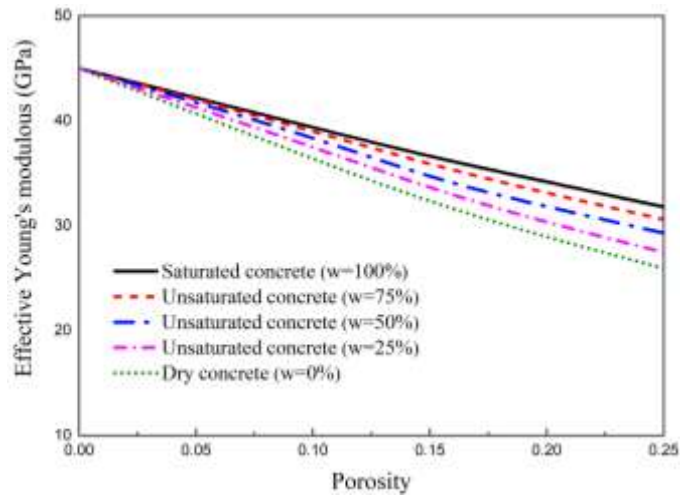

Fig. 17. Effective Young's modulus predicted by peridynamics with different saturations.

**4.2 Impact failure against a 300mm concrete target**

The impact perforation of saturated concrete targets is analyzed through the presented model, as experimentally tested in (Forquin et al., 2015). The same geometry and loading configuration from experiments are chosen in the simulations, as shown in Fig. 18. The target is in a steel culvert. The projectile is with a mass of 2.44kg and a velocity of 333m/s. In the numerical simulation, the projectile is assumed without damage during the impact process. To represent the stiff properties of the projectile, we adopted a particularly great value of Young's modulus here. The compressive strength of concrete is 39.5MPa, and the porosity is chosen to be 6%. Concrete parameters are given in Table 1 (Li et al., 1999). Other parameters include $P_{crush}$=14MPa, $\mu_{crush}$=8.1×$10^{-4}$, $P_{crush}$ = 3GPa, $\mu_{lock}$ = 0.16, $k_1$ = 15.7GPa, $k_2$ = -30.8GPa, and $k_3$ = 10.8GPa. For peridynamic modeling, the horizon size is 28mm, and we use a uniform discretization $\delta=4\Delta x$. The total number of nodes is 579,872.

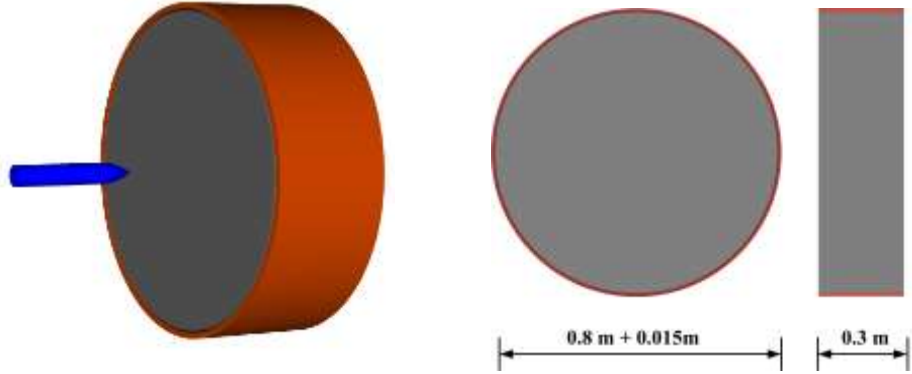

Fig. 18. The geometric model for the concrete impact perforation.

Table. 1 The material property of each phase in concrete.

| Concrete phase | Young's modulus (GPa) | Fracture energy (N/m) |
| --- | --- | --- |
| concrete | 32.0 | 107.0 |
| aggregate | 56.5 | 365.0 |
| cement mortar | 26.3 | 110.0 |
| ITZ | 20.2 | 90.0 |

The numerical results in terms of projectile acceleration and velocity are illustrated in Fig. 19, where the corresponding experimental data is shown as well. In Fig. 19(a), the projectile acceleration in numerical results is slightly higher than that of experimental data in the middle part, and the corresponding projectile velocity of the former is slightly smaller than that of the latter (shown in Fig. 19(b)). The apparent difference between numerical and experimental results is focused on the middle area. Moreover, the predicted residual velocity attained at 2ms is about 82m/s, which is a little greater than 75m/s measured from the experimental tests. Though there are some small divergences between numerical results and experimental data, a good consistency could still be observed, demonstrating the capability of the presented model and method for the impact failure of saturated concrete.

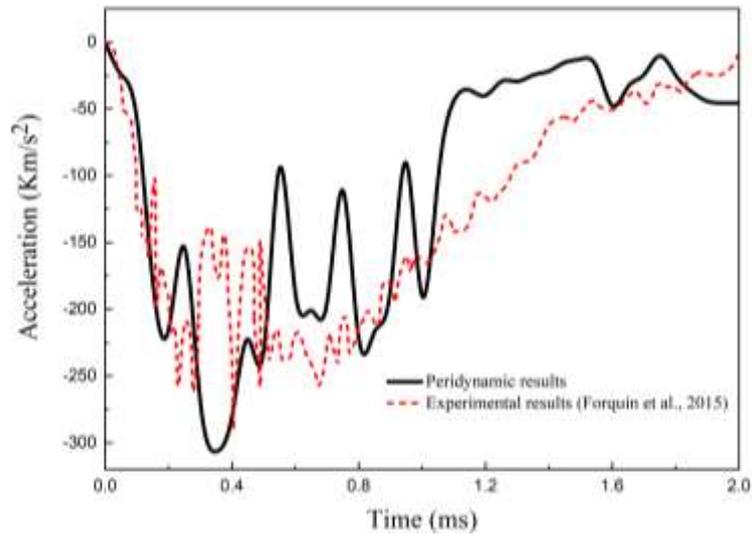

(a)

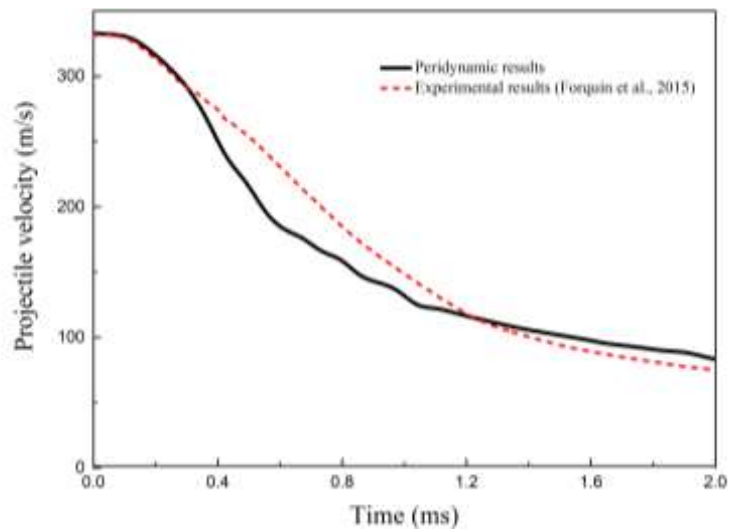

(b)

Fig. 19. Comparisons between experimental data and numerical predictions for a 300mm concrete target. (a) Projectile acceleration; (b) Projectile velocity.

Fig. 20 shows the comparison of damage patterns observed from the experimental tests (c) and peridynamic simulations. We can see that the damage pattern of the impact surface measured in the experiment is asymmetric, while the failure in the upper left corner is significantly different. This phenomenon may be caused by the small yaw angle of the projectile (Forquin et al., 2015). From Fig. 20(a) and Fig. 20(b), the failure modes predicted by the presented model are similar to those from the experimental observations. The crater radius and depth of the target observed experimentally are about 223.67mm and 84.91mm, and the scabbing radius and depth are about 375.26mm and 146.38mm, respectively. Similarly, the crater radius and depth of the target are predicted by peridynamics are about 257.43mm and 58.33mm, and the scabbing radius and depth are about 370.49mm and 141.04mm, respectively.

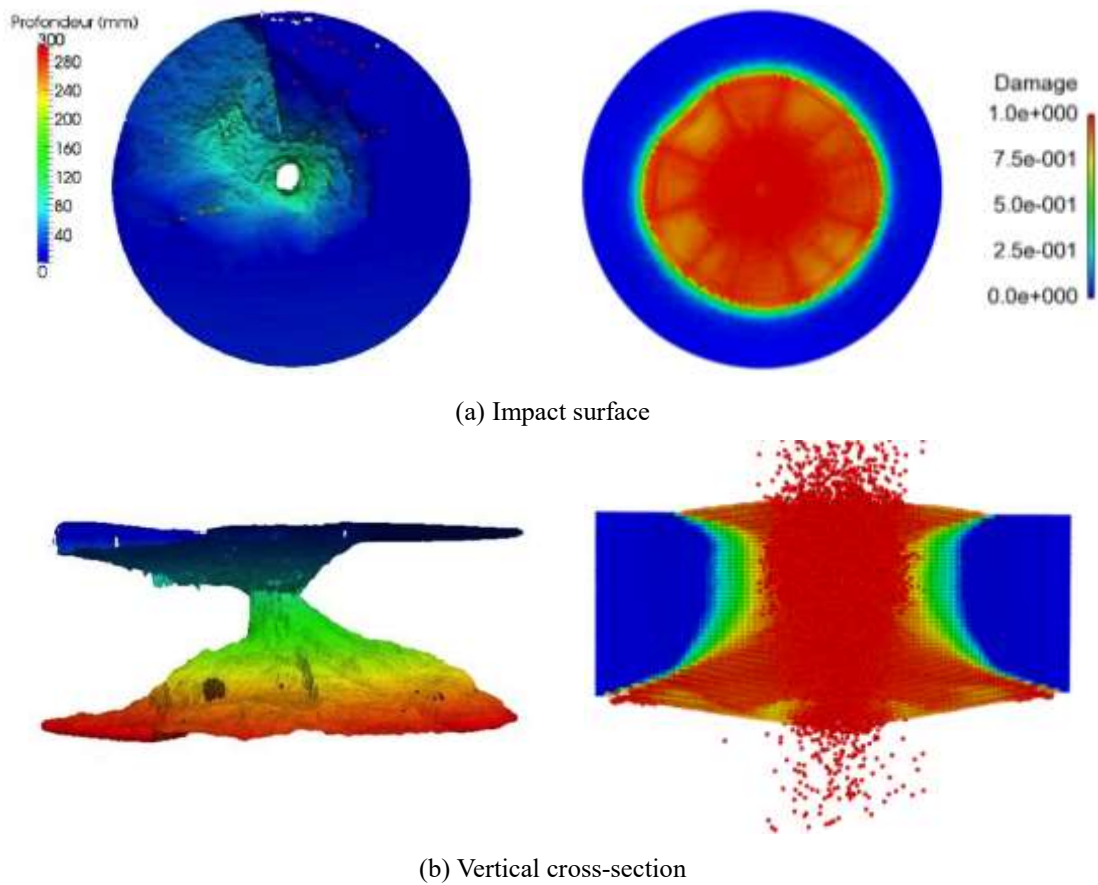

(a) Impact surface

(b) Vertical cross-section

Fig. 20. Damage maps of saturated concrete in experimental observations (left) and peridynamic predictions (right).

In the aforementioned example, the concrete is regarded as the saturated case. To discuss the influence of saturation on wet concrete, several other cases are conducted with different saturations for the impact perforation simulation.

Fig. 21 shows the influence of saturation on the projectile acceleration and residual velocity during the perforation. Three types of wet concrete are presented as saturated ($\varphi=1$), half-saturated ($\varphi=0.5$), and dry ($\varphi=0$). As shown in Fig. 21(a), the acceleration curves are basically identical at the beginning. The difference between these curves is much more evident after 0.5ms. A similar trend for the projectile residual velocity could be observed in Fig. 21(b). Moreover, a significant difference is the projectile residual velocity among different saturations, where the highest in dry concrete, the lowest in saturated concrete, and half-saturated concrete lie between. It can be concluded that the strain rate effect of saturated concrete is more sensitive than that of dry concrete, due to the inertia effect of the concrete itself and the viscous effect of internal free water, and results in a higher penetration resistance for wet concrete when subjected to the high-velocity impacting load.

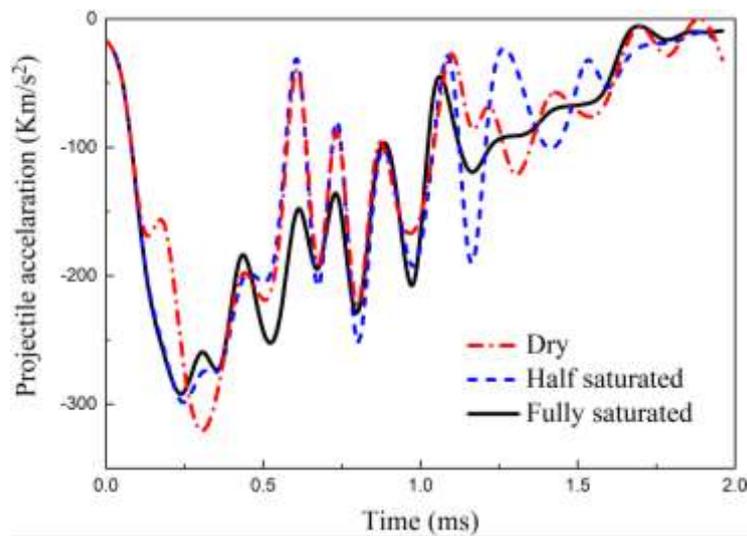

(a)

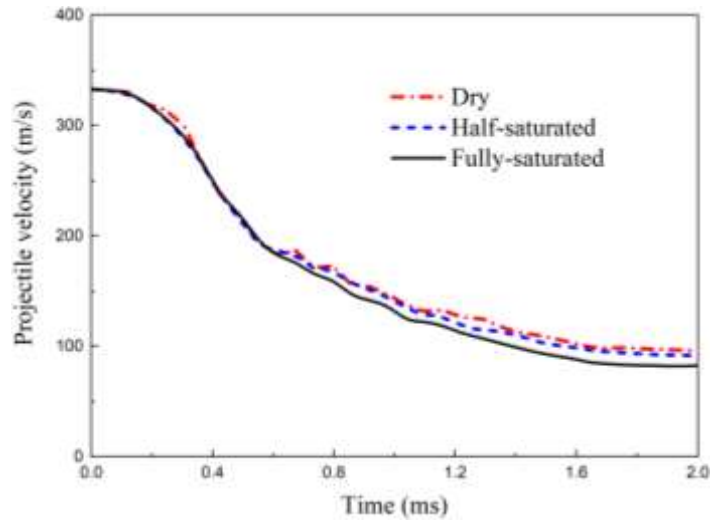

(b)

Fig. 21. Numerical results by considering saturated, half-saturated, and dry concrete. (a) Projectile acceleration; (b) Projectile velocity.

The failure modes of wet concrete associated with different saturations are illustrated in Fig. 22, which do not change significantly. Nevertheless, the cratering and scabbing sizes are a bit different, and the damaged area gradually enlarges when the saturation is higher. The crater radius, crater depth, and scabbing depth of impact surface for saturated concrete are 257.43mm, 58.33mm, and 141.04mm, the half-saturated concrete is 252.52mm, 56.44mm, and 111.93mm, and the dry concrete is 233.106mm, 45.16mm, and 98.5mm.

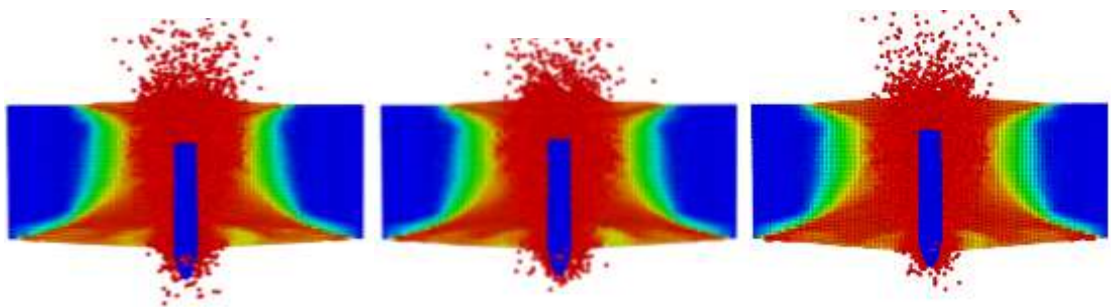

(a) Dry concrete  (b) Half-saturated concrete  (c) Saturated concrete

Fig. 22. The effect of wet concrete's saturation on crack patterns of the concrete target.

**4.3 Impact failure against an 800mm concrete target**

In this section, the geometric shape and boundary condition of the concrete target are the same as those conducted in the previous section, except for the thickness of 800mm, as shown in Fig. 23. The projectile is with a mass of 2.44kg and an initial velocity of 347m/s. In the peridynamic simulation, the horizon is $\delta=36\text{mm}$, which is four times the grid spacing. All of the material parameters are identical to the last section as well. There are 807,281 nodes in total.

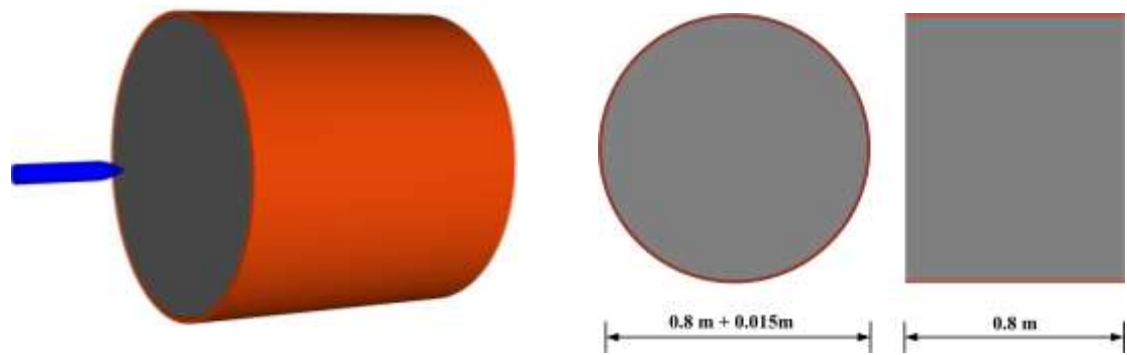

Fig. 23. The geometric model for the concrete impact perforation.

Numerical results and those measured experimentally in terms of the projectile acceleration and penetration depth are shown in Fig. 24. From Fig. 24(a), it could be observed that the projectile acceleration, predicted by the presented model, gradually reduces after 0.35ms. In the experimental measurements, the projectile acceleration rapidly drops to zero at about 1.4ms. Due to the difference in projectile acceleration, the variation trend of the penetration depth is different as well, especially after 0.35ms, as shown in Fig. 24(b). Nevertheless, the final penetration depth predicted by peridynamics is highly consistent with the experimental data, that the numerical result is about 0.261m, while the measured value is about 0.267m. Fig. 25 presents the comparison of damage patterns observed from the experimental measurements (Forquin et al., 2015) and numerical predictions. There is a satisfactory agreement between the numerical and experimental values for the crater radius. The average crater radius provided in the test is about 315.26mm, and the result predicted by the present model is about 323.92mm. These foregoing great agreements demonstrate again the applicability of the presented model for such impact problems of wet concrete.

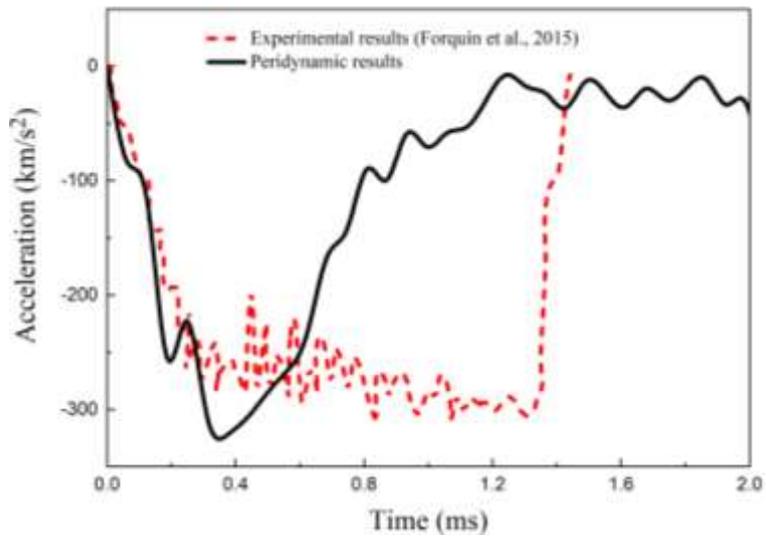

(a)

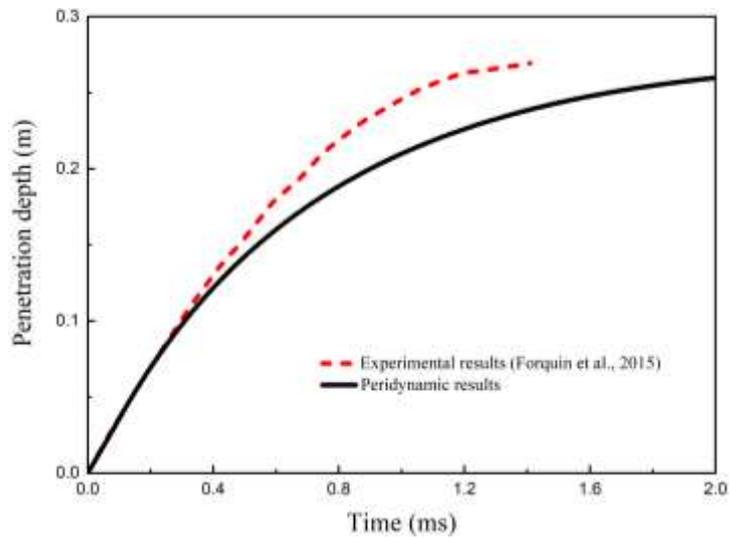

(b)

Fig. 24. Comparisons between experimental data and numerical predictions for the 800mm-thick wet concrete. (a) Projectile acceleration; (b) Projectile velocity.

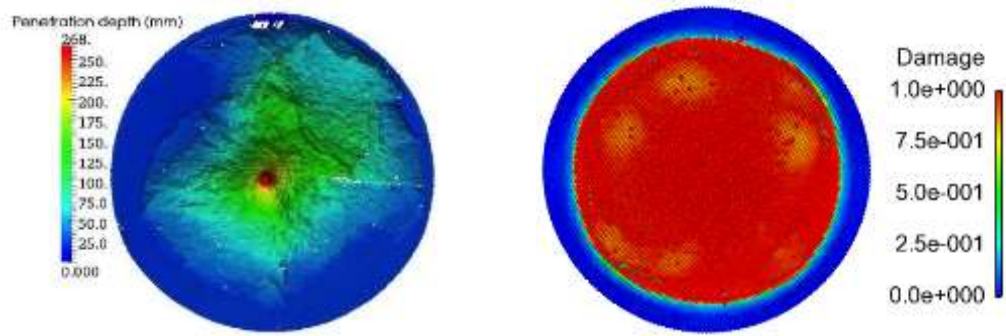

Fig. 25. Damage maps of the 800mm-thick saturated concrete in experimental observations (Forquin et al., 2015) (left) and peridynamic predictions (right).

The penetration depth in wet concrete with different saturations is shown in Fig. 26.

Similar to the previous section, the highest value of penetration depth is in dry concrete, while the lowest value is in saturated concrete, and half-saturated concrete lies between. We can conclude that wet concrete has a higher bearing capacity subjected to high-velocity impact, but a lower one for dry concrete.

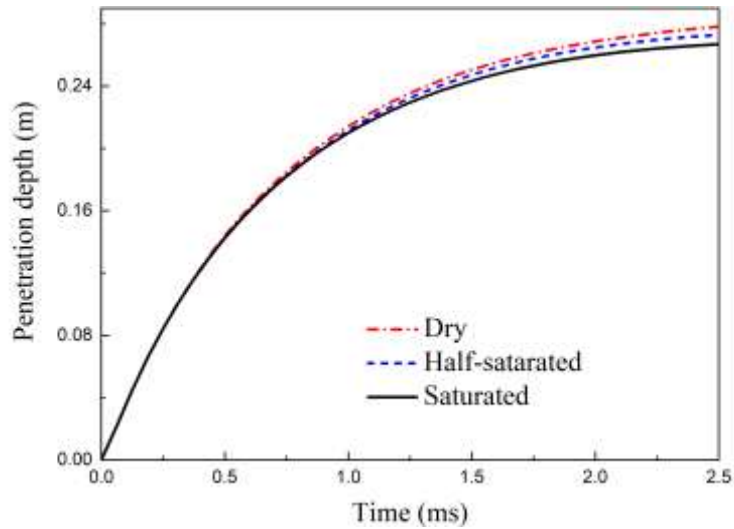

Fig. 26. Numerical results for penetration depth by considering saturated, half-saturated, and dry concrete.

The failure modes associated with different saturations are presented in Fig. 27. As can be seen from Fig. 27, the crater radius of the impact surface with different saturations is consistent with each other. However, the projectile position inside the target is inclined in saturated concrete, and there is an obvious cross-section damaged area as well. These phenomena are not observed in half-saturated and dry concrete cases.

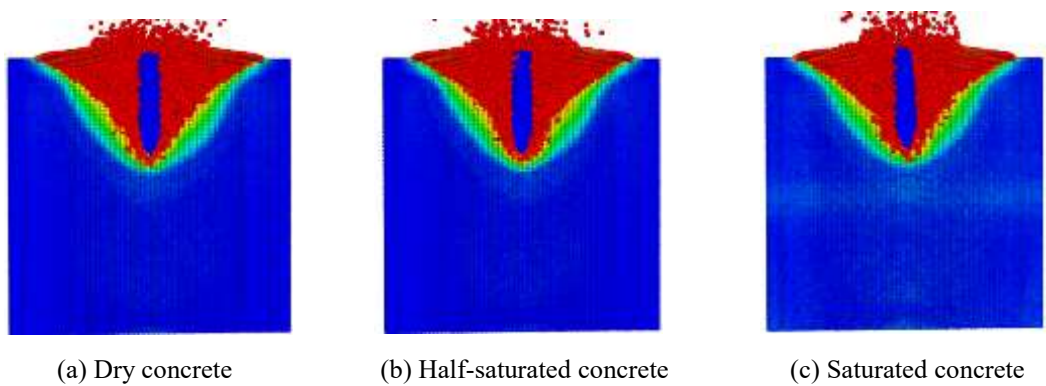

(a) Dry concrete  (b) Half-saturated concrete  (c) Saturated concrete

Fig. 27. The effect of wet concrete's saturation on crack patterns of the concrete target.

## 5. Conclusions

In this paper, a modified intermediately-homogenized peridynamic (IH-PD) model for wet concrete is established. This model considers the effect of porosity and free water on the mechanical behavior of wet concrete, where both saturated and unsaturated concrete is analyzed. The numerical simulation of effective Young's modulus and impact failure for wet concrete are conducted and show sensible agreements with corresponding experimental data:

(1) The modified IH-PD model presented in this paper can effectively describe the dynamic mechanical behavior of wet concrete, and characterize the influence of different saturations on the impact failure.

(2) The effective Young's modulus of wet concrete decreases with the increase of porosity, no matter what the saturation is. The free water inside the pores limits the deformation of the mortar matrix around, resulting in an improvement in the concrete stiffness.

(3) When subjected to the high-velocity impacting loads, the residual velocity of projectile in dry concrete is higher than that of wet concrete, and the penetration depth in dry concrete is also greater. It can be concluded that wet concrete displays a better penetration resistance for impact cases, due to the more sensitive characteristics to strain rate.

**Acknowledgments:**

The authors acknowledge the support of the National Natural Science Foundation of China (No.12072104, 51679077).